\documentclass[oneside,reqno,12pt]{amsart}
\usepackage[greek,english]{babel}
\usepackage{amsthm}
\usepackage{amsbsy}
\usepackage{amsfonts}
\usepackage{graphicx}
\usepackage{hyperref}
\hypersetup{colorlinks=true, citecolor=blue, linkcolor=red}

\newtheorem{thm}{Theorem}

\newtheorem{pro}{Proposition}
\newtheorem{rem}{Remark}

\numberwithin{equation}{section} \numberwithin{lem}{section}
\numberwithin{thm}{section} \numberwithin{cor}{section}
\numberwithin{pro}{section} \numberwithin{rem}{section}

\begin{document}

\title[Linearized inner problem for the phase separation of BEC's]{On a linear equation arising in the study of phase separation of Bose-Einstein condensates}
\author{Christos Sourdis}
\address{General Lyceum of Moires, Crete.
}
\email{sourdis@uoc.gr}

\date{\today}
\begin{abstract}
We consider the inner limit  system describing the phase separation in two-component Bose-Einstein condensates  linearized around the
one-dimensional solution
 in an infinite strip with zero and periodic boundary conditions, and obtain optimal invertibility estimates for the Fourier modes without necessarily assuming orthogonality conditions.
 \end{abstract}
 \maketitle

\section{Introduction}

The phase separation of a binary mixture of Bose-Einstein condensates is modelled by a Gross-Pitaevskii system and has been studied extensively when the
intercomponent repulsion is very strong (we refer to \cite{kowalczykPistoiaVaira} and the references therein). This parameter regime induces a singular perturbation problem. Solitary wave solutions satisfy a coupled semilinear
elliptic system, where segregated nodal domains divided by an interface appear in
the singular limit. The limiting components have disjoint supports in the aforementioned nodal domains. In the language of singular perturbation theory (see
for instance \cite{miller}) they are called \emph{outer solutions}. The regularity properties of the
outer solutions and the corresponding interface have been studied extensively in
\cite{noris,hugo}. Generically, for very strong repulsion, the behaviour across the interface
is governed by a suitable translated and scaled one-dimensional solution of (\ref{eqbereg})
below (see \cite{soaveZilioAHP}), the \emph{inner solution}.

Our motivation for the current paper is the study of the associated linearized operator at the aforementioned solutions, for large values of the intercomponent repulsion,
and its invertibility properties. This may have numerous applications ranging from studying the stability properties of such solutions to the construction of solutions with desired segregated asymptotic profiles for the singularly perturbed elliptic system.
For simplicity, let us assume that the interface is sufficiently smooth. Then, based on \cite{kowalczykPistoiaVaira} and accumulated evidence
from related problems, we have to first understand the linearized inner limit problem (the PDE version of (\ref{eqbereg})) on its one dimensional solution in an infinite  strip with zero and periodic boundary conditions.
  To this end, it is natural to perform a Fourier decomposition and arrive at the ODE problems (\ref{eqDir}) below with \[R=1/\epsilon\ \textrm{and}\ \omega=\epsilon \sqrt{\lambda_k},\]
  where $\epsilon>0$ is the small singular perturbation parameter (related to the large intercomponent repulsion) and $\lambda_k$ are the eigenvalues of a Schr\"{o}dinger operator on the interface (see \cite{kowalczykPistoiaVaira}).

We consider the linearised  operator
\[
L_\omega=\left( \begin{array}{cc}
                  -\frac{\partial^2}{\partial^2x}+V_2^2 & 2V_1V_2 \\
                  2V_1V_2 & -\frac{\partial^2}{\partial^2x}+V_1^2
                \end{array}
\right)+\omega^2 \textrm{Id},
\]
where $\omega\geq 0$ and $V_1,V_2>0$ satisfy
\begin{equation}\label{eqbereg}
\left\{
  \begin{array}{c}
    -V_1''+V_1V_2^2=0\ \textrm{in}\ \mathbb{R}, \\
    -V_2''+V_2V_1^2=0\ \textrm{in}\ \mathbb{R}, \\
  \end{array}
\right.
\end{equation}
$V_1(0)=V_2(0)=1$, $V_1(-x)=V_2(x)$, $V_1'(x)>0$, $x\in \mathbb{R}$, and
\begin{equation}\label{eqasym}\begin{array}{ll}
                                V_1(x)=Ax+B+\mathcal{O}(e^{-cx^2}) & \textrm{as}\ x\to +\infty, \\
                                V_1(x)=\mathcal{O}(e^{-cx^2}) & \textrm{as}\ x\to -\infty,
                              \end{array}
\end{equation}
for some constants $A>0$, $B\in \mathbb{R}$ and $c>0$
(the existence and uniqueness of such a $(V_1,V_2)$ has been established in \cite{berestArma} and \cite{berestAdv}, respectively).

The following Dirichlet problem was studied in \cite[Sec. 2]{kowalczykPistoiaVaira}:
\begin{equation}\label{eqDir}\begin{array}{rc}
                               L_\omega \varphi=g, &  \ \textrm{in}\ (-R,R), \\
                               \varphi(\pm R)=0, &
                             \end{array}
\end{equation}
where $R>0$. In Subsection 2.3 of the aforementioned reference, \emph{a uniform estimate in $\omega \in [0,\infty)$} was obtained for (\ref{eqDir}) which, in particular, implies that
\begin{equation}\label{eqprop}
|\varphi(0)|\to 0  \ \textrm{if}\ \|g\|_{C^0_\theta\left((-R,R)\right)}\to 0\ \textrm{as}\ R\to \infty,
\end{equation}
where
\begin{equation}\label{eqnorm}
\|g\|_{C^0_\theta\left((-R,R)\right)}=\|g \cosh(\theta x)\|_{C^0\left((-R,R)\right)},
\end{equation}
and $\theta>0$ (independent of $\omega, R$). This was accomplished by constructing a suitable positive supersolution to (\ref{eqDir}), i.e., a pair $\bar{u}(x)=\left(\bar{u}_1(x),\bar{u}_2(x)\right)\geq 0$ such that
\[
  L_\omega \bar{u}\geq g,\ (\textrm{componentwise});
\]
we stress that $\bar{u}$ was chosen independently of  $\omega$. Then, the authors appealed to a maximum principle type result of \cite{figu} concerning \emph{cooperative} linear elliptic systems. However, according to our understanding, the system at hand becomes cooperative when written in terms of $(\varphi_1,-\varphi_2)$
(see also \cite[Sec. 2]{berestAdv}).
Indeed, after this transformation, the off diagonal terms in the resulting righthand sides of (\ref{eqDir}) have negative coefficients, as required by \cite{figu} in order for the system to be considered  cooperative. In any case, this was only needed in \cite{kowalczykPistoiaVaira} in order to give a different proof of a result of \cite{aftalionSour} on the solvability of
\begin{equation}\label{eqentire}
  L_0\varphi=g\  \textrm{in}\  \mathbb{R}\ \textrm{with}\ \|g\|_{C_\theta^0(\mathbb{R})}<\infty,
  \end{equation}
  with the following orthogonality conditions
  \begin{equation}\label{eqorthdouble}
    \int_{-\infty}^{+\infty}(V_1'g_1+V_2'g_2)dx=0\ \textrm{and}\ \int_{-\infty}^{+\infty}\left((xV_1'+V_1)g_1+(xV_2'+V_2)g_2\right)dx=0,
  \end{equation}
  where $(V_1',V_2')$, $(xV_1'+V_1,xV_2'+V_2)$ belong in the kernel of $L_0$;
which subsequently played an important role in \cite{kowalczykPistoiaVaira} for the use of the \emph{infinite-dimensional Lyapunov-Schmidt reduction} and a gluing argument in a class of semilinear elliptic systems modeling phase separation in two-component Bose-Einstein condensates.

The purpose of the current paper is twofold. Firstly, we will provide a simple counterexample to (\ref{eqprop}) for $0<\omega\ll 1$ depending on $R\gg 1$ (this is the most interesting regime since $\sigma(L_\omega)\subseteq [\omega^2,\infty)$, see \cite{berestArma,kowalczykPistoiaVaira}). This will be done in Section \ref{seccounter}. The starting idea is that, in light of (\ref{eqasym}) and a standard barrier argument applied to $\varphi_i$ for $(-1)^ix>0$, $i=1,2$, for $R\gg 1$ and $0<\omega\ll 1$, problem  (\ref{eqDir}) with $g$ negligible decomposes in two problems:
\begin{itemize}
  \item \underline{the inner one}
  \begin{equation}\label{eqpb1}
    \left\{\begin{array}{c}
             -\varphi_1''+V_2^2\varphi_1+2V_1V_2\varphi_2=0, \\
             -\varphi_2''+V_1^2\varphi_2+2V_1V_2\varphi_1=0,
           \end{array}
     \right.
      \end{equation}
      which provides an effective approximation to (\ref{eqDir}) for $|x|< M$ with $M$ large but independent of $R,\omega$;
  \item \underline{the outer ones}
\begin{equation}\label{eqpb2}
  -\varphi_1''+\omega^2\varphi_1=0,\ \varphi_1(R)=0,\ \varphi_2=0,\ 1 \ll x <R,
\end{equation}
and the corresponding one for $-R<x\ll -1$.
\end{itemize}
The inner problem (\ref{eqpb1}) admits a bounded solution $(V_1',V_2')$ which, however, does not satisfy the Dirichlet boundary conditions in (\ref{eqDir}). On the other hand, the outer solutions do satisfy the desired boundary conditions and, as it turns out, can be matched sufficiently well with the aforementioned inner solution $(V_1', V_2')$, provided that $\omega\ll 1$, yielding an element in the approximate kernel of $L_\omega$.

In the second and main part of the paper, having the construction of our counterexample in mind, we will establish  new a-priori estimates for (\ref{eqDir}).
To this end, we will use \cite{aftalionSour} as a guideline   which dealt with the case $\omega=0$. However, some nontrivial modifications are needed, the main one being on
how to connect (exchange) the information between the inner and outer regions (the latter are actually chosen strikingly different from \cite{aftalionSour}). We note in passing that this is the objective of the so called \emph{exchange lemmas} from {geometric singular perturbation theory} (see \cite{kuehn} and the references therein). However, we have not attempted to make a formal connection with the aforementioned theory. Nevertheless, here we employed a problem specific argument based on the self-adjoint property of $L_\omega$. Interestingly enough and in contrast to \cite{aftalionSour}, our proof does not involve the use of the other element of the kernel of $L_0$ that appears in the second relation of (\ref{eqorthdouble}).

 We insist  that establishing a-priori estimates for (\ref{eqDir}) or (\ref{eqentire}), with or without the orthogonality conditions (\ref{eqorthdouble}), is of great importance in the study of the corresponding linear PDE problems in a strip. In fact, the former ODE problems  arise after performing a  Fourier
decomposition in the latter PDE ones (see \cite{kowalczykPistoiaVaira}). The authors of \cite{kowalczykPistoiaVaira} introduced the orthogonality conditions in (\ref{eqorthdouble}) in order to deal with the
problem of \emph{small divisors} of the aforementioned linearized  inner PDE (employing the corresponding result of \cite{aftalionSour}  for the case $\omega=0$). In this regard, we emphasize that our main result  below estimates precisely those small divisors in terms of $0<\omega\ll 1$.

Our main result is the following.
\begin{thm}\label{thm}
Given $\theta>0$, there exist constants $C, R_0, \omega_0, P_0>0$ such that
\begin{equation}\label{eqDirThm}
L_\omega \varphi=g,\ x\in (-R,R),\ \varphi(\pm R)=0,
\end{equation}
with
\[R>R_0,\ \omega\in (0,\omega_0) \ \textrm{and}\ \omega R>P_0\]
implies that
\begin{equation}\label{eqthm}
  \|\varphi\|_{C^0\left((-R,R)\right)}\leq C \omega^{-1}\|g\|_{C^0_\theta\left((-R,R)\right)},
\end{equation}
where the norm in the righthand side of the above estimate was defined in (\ref{eqnorm}).
\end{thm}
It is worth mentioning that the counterexample for (\ref{eqprop}) that we will construct in Section \ref{seccounter}, see Proposition \ref{procounter} therein, indicates that (\ref{eqthm}) is  sharp.

We point out that if $\omega R\leq P_0$, then   Theorem \ref{thm} is a direct consequence of \cite[Sec. 4]{casteras}.
Moreover, as can be apparent from the proofs, all the above continue to hold if the problem is posed in the entire real line.

For completeness purposes and for comparison, let us mention that if $\omega=0$ it was shown in \cite[Prop. 2.2]{aftalionSour} that the corresponding assertion of Theorem
\ref{thm} is
\[
  \|\varphi\|_{C^0\left((-R,R)\right)}\leq C R\|g\|_{C^0_\theta\left((-R,R)\right)}.
\]
Furthermore, if additionally one of the orthogonality conditions in (\ref{eqorthdouble}) applies (with the integrals taken over $(-R,R)$), then
the stronger estimate
\[
  \|\varphi\|_{C^0\left((-R,R)\right)}\leq C \|g\|_{C^0_\theta\left((-R,R)\right)}
\]is valid.
In fact, at least in the case of the whole real line and still for $\omega=0$, when both of the orthogonality conditions in (\ref{eqorthdouble}) hold, we conclude that
\begin{equation}\label{eqkalergi}
  \|\varphi\|_{C^0_\theta\left(\mathbb{R}\right)}\leq C \|g\|_{C^0_\theta\left(\mathbb{R}\right)},
\end{equation}
(see \cite[Prop. 2.3]{aftalionSour} and \cite[Lem. 2.2]{kowalczykPistoiaVaira}).

In \cite[Sec. 4]{kpwGelfand}, the authors studied the solvability of a scalar linearized  problem involving an exponential nonlinearity.
On the one hand, this   shares  several common
features with the problem herein.  On the other hand, the corresponding $\omega$-problem to (\ref{eqDir})  therein admits  an explicit change of variables  which transforms the
  homogeneous equation   into one of the   Legendre family. Actually, their emphasis   was placed on righthand sides that
satisfy both of the analogous orthogonality conditions  to (\ref{eqorthdouble}).  For the problem at hand, however, since $(V_1,V_2)$ is not given
explicitly, we believe that a suitable transformation  is unlikely to exist.
Nevertheless, based on the    results of \cite{kpwGelfand}  and Remark \ref{remorth} below, we expect that property (\ref{eqkalergi})
continues to hold at least for small $\omega>0$.

The outline of the paper is the following. In Section \ref{seccounter}, we will provide the construction of our aforementioned counterexample. Lastly, in Section \ref{secmain}
we will prove our main result.

\section{Construction of an element in the approximate kernel of $L_\omega$: a counterexample to (\ref{eqprop})}\label{seccounter}
For $x\in [0,R]$, we let
\begin{eqnarray}\label{eqcounterDef}
   \nonumber
    \varphi_1(x) &=& \left(V_1'(x)-A\right)\eta_R(x)+A\frac{\sinh\left(\omega(R-x) \right)}{\sinh(\omega R)}, \\
    \varphi_2(x) &=& V_2'(x)\eta_R(x),
  \end{eqnarray}
where $A$ as in (\ref{eqasym}) and
\[
\eta_R(x)=\zeta\left(\frac{|x|}{\ln R} \right)
\]
with $\zeta$ a smooth cutoff function such that
\[
\zeta(y)=\left\{\begin{array}{ll}
                  1 & \textrm{if}\ |y|<\frac{1}{2}, \\
                  &\\
                  0 & \textrm{if}\ |y|>\frac{3}{4}.
                \end{array}
\right.
\] We note in passing that in (\ref{eqcounterDef}) we have employed a standard matching procedure from asymptotic analysis, see for instance \cite[Sec. 8.5]{miller} and recall the discussion related to (\ref{eqpb1}) and (\ref{eqpb2}).
For $x\in [-R,0]$, we let
\begin{equation}\label{eqsymmetry4}
\left(\varphi_1(x),\varphi_2(x) \right)=\left(-\varphi_2(-x),-\varphi_1(-x) \right).
\end{equation}

In the remainder of the proof, we will work with the assumption that
\begin{equation}\label{eqass}
  \omega R\geq 1.
\end{equation}
Moreover, throughout the paper, we will denote by $c/C$ a small/large generic constant which is independent of both $0<\omega\ll 1$ and $R\gg 1$, whose value may decrease/increase  as the proof progresses (the usual Landau symbols will be understood in the same sense).

For $x\in [0,(\ln R)/2]$, we have
\[
\varphi_1(x)=V_1'(x)+\sigma(x),\ \ \varphi_2(x)=V_2'(x),
\]
where
\begin{equation}\label{eqsigmaDef}
  \sigma(x)=A\left[\frac{\sinh\left(\omega(R-x) \right)}{\sinh(\omega R)}-1\right].
\end{equation}
It is easy to check that $\sigma(0)=0$ and $|\sigma'|\leq C\omega$ (recall (\ref{eqass})), which imply that
\begin{equation}\label{eqsigmaEstima}
  |\sigma(x)|\leq C\omega x,\ x\in \left[0, \frac{\ln R}{2} \right].
\end{equation}We also note that
\begin{equation}\label{eqsigmaEq}
-\sigma''+\omega^2\sigma =-A\omega^2.
\end{equation}
Hence, keeping in mind (\ref{eqasym}), we find
\[
L_\omega\left(\begin{array}{c}
                \varphi_1 \\
                \varphi_2
              \end{array}
 \right)=\left(\begin{array}{c}
                \omega^2V_1'-A\omega^2+\mathcal{O}(e^{-cx^2})\sigma \\
                \omega^2V_2'+\mathcal{O}(e^{-cx^2})\sigma
              \end{array}
 \right)\stackrel{(\ref{eqsigmaEstima})}{=}\mathcal{O}(\omega^2)+\omega \mathcal{O}(e^{-cx^2}),\ x\in\left[0,\frac{\ln R}{2} \right].
\]

For the remainder of the proof we will make the choice
\begin{equation}\label{eqomegachoice}
\omega=R^{-\alpha}\ \textrm{with}\ \alpha>0 \ \textrm{to be determined}.
\end{equation}

Then, given $\theta\geq 0$ (independent of $R$, $\alpha$), we find
\begin{equation}\label{eqcombo1}
\|L_\omega(\varphi_1,\varphi_2)\|_{C^0_\theta\left(0,\frac{\ln R}{2}\right)}\leq CR^{-2\alpha}R^{\frac{\theta}{2}}+CR^{-\alpha}.
\end{equation}

For $x\in \left[\frac{\ln R}{2},\frac{3\ln R}{4}\right]$, by virtue of (\ref{eqasym}) and (\ref{eqcounterDef}), we get
\begin{eqnarray}
   \nonumber
    \varphi_1(x) &=& A\frac{\sinh\left(\omega(R-x) \right)}{\sinh(\omega R)}+\mathcal{O}\left(e^{-c(\ln R)^2}\right), \\
    \varphi_2(x) &=& \mathcal{O}\left(e^{-c(\ln R)^2}\right).
  \end{eqnarray}
  We point out that the above relations hold uniformly with respect to $x$ and can be differentiated in the obvious manner. So, using (\ref{eqasym}) once more, we obtain
  \begin{equation}\label{eqcombo2}
  L_\omega\left(\begin{array}{c}
                  \varphi_1 \\
                  \varphi_2
                \end{array}
   \right)=\mathcal{O}\left(e^{-c(\ln R)^2}\right),\ \textrm{uniformly on}\ \left[\frac{\ln R}{2},\frac{3\ln R}{4}\right],\ \textrm{as}\ R\to \infty.
  \end{equation}
Lastly, since
\[
\varphi_1(x) = A\frac{\sinh\left(\omega(R-x) \right)}{\sinh(\omega R)},\
\varphi_2(x)=0,\ x\in
\left(\frac{3\ln R}{4},R\right),
\] we see that
\begin{equation}\label{eqcombo3}
L_\omega\left(\begin{array}{c}
\varphi_1 \\
\varphi_2
\end{array}
\right)=
\mathcal{O}\left(e^{-c(\ln R)^2}\right),\ \textrm{uniformly on}\
\ x\in \left[\frac{3\ln R}{4},R\right],\ \textrm{as}\ R\to \infty.
\end{equation}

By combining (\ref{eqcombo1}), (\ref{eqcombo2}), (\ref{eqcombo3}), and making the choice
\[
\alpha=\theta<1\ \textrm{in}\ (\ref{eqomegachoice}),
\]
which respects (\ref{eqass}), we deduce that
\[
\|L_\omega(\varphi_1,\varphi_2)\|_{C_\theta^0\left((0,R)\right)}\leq C R^{-\theta}.
\]
Clearly, thanks to the symmetry (\ref{eqsymmetry4}), the above estimate continues to hold over $(-R,0)$.

Consequently, we have arrived at the following proposition which provides a counterexample to the property (\ref{eqprop}).
\begin{pro}\label{procounter}
Given $\theta\in (0,1)$, the pair $(\varphi_1,\varphi_2)$ as defined in (\ref{eqcounterDef}) and (\ref{eqsymmetry4}) satisfies
\[
\limsup_{R\to \infty}\omega^{-1}\|L_\omega(\varphi_1,\varphi_2)\|_{C_\theta^0\left((-R,R)\right)}<\infty \ \textrm{with}\ \omega=R^{-\theta},
\]
while
\[
(\varphi_1,\varphi_2)\to(V'_1,V'_2)\ \textrm{in}\ C_{loc}(\mathbb{R})\ \textrm{as}\ R\to \infty.
\]
\end{pro}

\section{Proof of the main result}\label{secmain}

\subsection{Proof of Theorem \ref{thm}}
\begin{proof}
We will argue by contradiction. So, let us suppose that there exist sequences  $\varphi_n=(\varphi_{1,n},\varphi_{2,n})$, $g_n=(g_{1,n},g_{2,n})$ such that
\begin{equation}\label{eqproofDirThm}
L_{\omega_n} \varphi_n=g_n,\ x\in (-R_n,R_n),\ \varphi_n(\pm R_n)=0,
\end{equation}
with
\begin{equation}\label{eqproofomegaRa}
R_n\to \infty,\ \omega_n\to 0,\ \omega_n R_n\to \infty,
\end{equation}
\begin{equation}\label{eqproofContra}
  \|\varphi_n\|_{C^0\left((-R_n,R_n)\right)}=1\ \textrm{and}\  \omega_n^{-1}\|g_n\|_{C^0_\theta\left((-R_n,R_n)\right)}\to 0.
\end{equation}

Throughout the remainder of the proof, we will denote by $c/C$ a small/large generic
constant which is independent of large $n$, whose value
may decrease/increase as the proof progresses (the usual Landau symbols will be
understood in the same sense). Moreover, for notational convenience, we will frequently drop the subscript $n$.

\underline{The inner zone $|x|<M$ with $M\gg 1$ independent of $n\gg1$.}

 By (\ref{eqproofDirThm})-(\ref{eqproofContra}), standard elliptic estimates and the usual diagonal-compactness argument, passing to a subsequence if needed, we have
\begin{equation}\label{eqinnerloc1}
\varphi_n\to \varphi_\infty \  \textrm{in}\ C^1_{loc}(\mathbb{R}),
\end{equation}
where $\varphi_\infty$ satisfies
\[L_0\varphi_\infty=0\ \textrm{in}\ \mathbb{R}\ \textrm{(in the classical sense)\ and}\ \|\varphi_\infty\|_{C^0(\mathbb{R})}\leq 1.
\]
Then, by the nondegeneracy of $V=(V_1,V_2)$ (see \cite[Thm. 1.3]{berestArma}), we infer that
\begin{equation}\label{eqinnerMainLokatzis}
  \varphi_\infty\equiv \lambda V' \ \textrm{for some}\ \lambda\in \mathbb{R}.
\end{equation}

\underline{The outer zone $\omega_n^{-\beta}<|x|<R_n$ with $0<\beta\ll 1$ independent of $n\gg1$.}

Let us first emphasize that our outer zone    differs considerably from that in the case $\omega=0$, see \cite[Prop. 2.2]{aftalionSour}, where it was plainly chosen as $\delta R<|x|<R$ with $\delta$ small but independent of $n$.

Now, let \[\beta \in (0,1)\]   be arbitrary but independent  of $n\gg1$. We note that (\ref{eqproofomegaRa}) implies that
\begin{equation}\label{eqfrac}
\frac{\omega_n^{-\beta}}{R_n}\to 0.
\end{equation}
If $\frac{1}{2}\omega_n^{-\beta}<x<R_n$, then (\ref{eqasym}), (\ref{eqproofDirThm}) and (\ref{eqproofContra}) yield
\begin{equation}\label{eqoutEq}
-\varphi_1''+\omega^2\varphi_1=E\ \textrm{with}\ E=\mathcal{O}\left(e^{-\frac{c}{\omega^\beta}} \right)\ \textrm{uniformly in}\ x;\ \varphi_1(R)=0.
\end{equation}
Let $\hat{\varphi}_1$ be uniquely determined from
\begin{equation}\label{eqhatosEq}
-\hat{\varphi}_1''+\omega^2\hat{\varphi}_1=E\ \textrm{in}\ \left(\frac{1}{2}\omega^{-\beta},R \right);\ \hat{\varphi}_1\left(\frac{1}{2}\omega^{-\beta}\right)=\hat{\varphi}_1(R)=0.
\end{equation}
By the maximum principle and the bound of $E$ from (\ref{eqoutEq}),  we deduce that
\begin{equation}\label{eqhatunifW}
\|\hat{\varphi}_1\|_{C^0\left(\frac{1}{2}\omega^{-\beta},R \right)}\leq \frac{1}{\omega^2}\|E\|_{C^0\left(\frac{1}{2}\omega^{-\beta},R \right)}= \mathcal{O}\left(e^{-\frac{c}{\omega^\beta}}\right).
\end{equation}
In turn, by standard elliptic regularity estimates up to $x=R$ (applied to intervals of size $1$), we get
\begin{equation}\label{eqhatGradunifW}
\|\hat{\varphi}'_1\|_{C^0\left(\omega^{-\beta},R \right)}= \mathcal{O}\left(e^{-\frac{c}{\omega^\beta}}\right).
\end{equation}

Next, we let
\begin{equation}\label{eqtilda}
  \tilde{\varphi}_1=\varphi_1-\hat{\varphi}_1,\ x\in \left(\omega^{-\beta},R \right).
\end{equation}
From (\ref{eqoutEq}) and (\ref{eqhatosEq}) we obtain
\begin{equation}\label{eqtildaEq74}
-\tilde{\varphi}_1''+\omega^2\tilde{\varphi}_1=0\ \textrm{in}\ \left(\omega^{-\beta},R \right);\ \tilde{\varphi}_1(R)=0.
\end{equation}
It follows that
\begin{equation}\label{eqtildemuSinh}
\tilde{\varphi}_1(x) = \mu_n\frac{\sinh\left(\omega(R-x) \right)}{\sinh\left(\omega (R-\omega^{-\beta})\right)}
,\ x\in (\omega^{-\beta},R),\ \textrm{with}\ |\mu_n|\leq 1.
\end{equation}
So, by combining (\ref{eqhatunifW}), (\ref{eqhatGradunifW}) and (\ref{eqtilda}), we find
\begin{equation}\label{eqtildemuSinhplus}
{\varphi}_1(x) = \mu_n\frac{\sinh\left(\omega(R-x) \right)}{\sinh\left(\omega (R-\omega^{-\beta})\right)}+\mathcal{O}_{C^1}\left(e^{-\frac{c}{\omega^\beta}}\right)
,\ x\in (\omega^{-\beta},R),\ \textrm{with}\ |\mu_n|\leq 1.
\end{equation}
Analogously, we can show that
\begin{equation}\label{eqtildenuSinhminus}
\varphi_2(x) = \nu_n\frac{\sinh\left(\omega(R+x) \right)}{\sinh\left(\omega (R-\omega^{-\beta})\right)}+\mathcal{O}_{C^1}\left(e^{-\frac{c}{\omega^\beta}}\right)
,\ x\in (-R,-\omega^{-\beta}),\ \textrm{with}\ |\nu_n|\leq 1.
\end{equation}
For future reference, we note that
\begin{equation}\label{future}
\varphi_i'\left((-1)^{i-1} \omega^{-\beta}\right) = \kappa_i\omega\frac{\cosh\left(\omega(R- \omega^{-\beta}) \right)}{\sinh\left(\omega (R-\omega^{-\beta})\right)}+\mathcal{O}\left(e^{-\frac{c}{\omega^\beta}}\right)\stackrel{(\ref{eqfrac})}{=} \mathcal{O}(\omega),\ i=1,2,
\end{equation}
where $\kappa_1=\mu_n$ and $\kappa_2=\nu_n$.
In fact, we have
\begin{equation}\label{futureprecise}
\varphi_2'(- \omega^{-\beta}) = \nu_n\omega+o(1)\omega\ \textrm{and}\ \varphi_1'( \omega^{-\beta}) =- \mu_n\omega+o(1)\omega\ \textrm{as}\ n\to \infty.
\end{equation}

\underline{The intermediate zone $M< |x|<\omega_n^{-\beta}$.} Our next objective is to bridge the gap $M<|x|<\omega^{-\beta}$ between the inner and outer zones. For the case $\omega=0$, this was accomplished in \cite[Prop. 2.2]{aftalionSour} by starting from the outer zone and working towards the inner one (solving the corresponding initial value problem for each component separately). In contrast, here we will essentially take the opposite route,  after testing the ODE system (\ref{eqproofDirThm}) with $(V_1',V_2')$ and  integrating by parts, which yields a scalar first order ODE (information from the outer region will be needed for some boundary terms). We believe that here lies  the main novelty in our proof.

Motivated from (\ref{eqinnerMainLokatzis}), we write
\begin{equation}\label{eqlokatzisPsi}
\varphi_i=\lambda V_i'+\psi_i\ \textrm{with}\ \psi_i\to 0\ \textrm{in}\ C^{1}_{loc}(\mathbb{R})\ \textrm{as}\ n\to \infty,\ i=1,2.
\end{equation}

Then, testing (\ref{eqproofDirThm}) by $(V_1',V_2')$ over $(x,\omega^{-\beta})$, $x>0$, and integrating by parts, recalling that $L_0V'=0$, yields
\begin{equation}\label{eqagni1}\begin{split}
\sum_{i=1}^{2}\left(\psi_i'V_i'(x)-\psi_iV_i''(x)\right)-\sum_{i=1}^{2}\left(\psi_i'V_i'(\omega^{-\beta})-\psi_iV_i''(\omega^{-\beta})\right) &  \\
+\omega^2\sum_{i=1}^{2}\int_{x}^{\omega^{-\beta}}\varphi_iV_i'(y)dy                                    & =\sum_{i=1}^{2}\int_{x}^{\omega^{-\beta}}g_iV_i'(y)dy.
                               \end{split}
\end{equation}
We point out that we tacitly used that\[
\varphi_i'V_i'(x)-\varphi_iV_i''(x)=\psi_i'V_i'(x)-\psi_iV_i''(x).
\]

In what follows we will estimate the various terms of the above relation. Firstly, by (\ref{eqasym}) and (\ref{eqlokatzisPsi}) we get
\begin{equation}\label{eqagni22}\psi_2'V_2'(x)-\psi_2V_2''(x)=o(1)e^{-cx^2},\ \textrm{uniformly on} \ [0,R_n],\ \textrm{as}\ n\to \infty.
\end{equation}
Next, thanks to (\ref{eqasym}), (\ref{future}) and (\ref{eqlokatzisPsi}), we find
\begin{equation}\label{eqagni3}
\psi_1'V_1'(\omega^{-\beta})=\mathcal{O}(\omega)\ \textrm{as}\ n\to \infty.
\end{equation}
Furthermore, from (\ref{eqasym}), (\ref{eqproofContra}) and (\ref{eqlokatzisPsi}) we obtain
\begin{equation}\label{eqagni4}
\psi_iV_i''(\omega^{-\beta})=\mathcal{O}\left(e^{-\frac{c}{\omega^\beta}}\right),\ i=1,2.\end{equation}
Similarly, since \begin{equation}\label{eqgradg}
                   |\varphi_i'|\leq C\ \textrm{on}\ [-R,R],\ i=1,2,
                 \end{equation} (from (\ref{eqproofDirThm})-(\ref{eqproofContra}) and standard elliptic estimates), recalling (\ref{eqlokatzisPsi}), we have
\begin{equation}\label{eqagni5}
\psi_2'V_2'(\omega^{-\beta})=\mathcal{O}\left(e^{-\frac{c}{\omega^\beta}}\right).
\end{equation}
Moreover, owing to (\ref{eqasym}) and (\ref{eqproofContra}), we can easily estimate the  integrals in (\ref{eqagni1}) in the following way:
\begin{equation}\label{eqagni6}
\omega^2\int_{x}^{\omega^{-\beta}}\varphi_iV_i'(y)dy=\mathcal{O}(\omega^{2-\beta})\ \textrm{and}\ \int_{x}^{\omega^{-\beta}}g_iV_i'(y)dy=o(\omega)e^{-\theta x},\ x\in (0,\omega^{-\beta}).
\end{equation}

Now, relation (\ref{eqagni1}) via (\ref{eqagni22})-(\ref{eqagni6}) yields
\begin{equation}\label{eqagni7}
\psi_1'V_1'(x)-\psi_1V_1''(x)=o(1)e^{-\theta x}+\mathcal{O}(\omega)+\mathcal{O}(\omega^{2-\beta}),\ \textrm{uniformly on}\ [0,\omega^{-\beta}],\ \textrm{as}\ n\to \infty.
\end{equation}
Since $V_1'(x)\geq V_1'(0)\geq c$ for $x\geq 0$, the above relation gives
\begin{equation}\label{eqagni8}
Q'(x)=\frac{\psi_1'V_1'-\psi_1V_1''}{\left(V_1'\right)^2}=o(1)e^{-\theta x}+\mathcal{O}(\omega)+\mathcal{O}(\omega^{2-\beta}),\ \textrm{where}\ Q=\frac{\psi_1}{V_1'},\ x\in (0,\omega^{-\beta}).
\end{equation}
Integrating the above relation in $(0,x)$, keeping in mind that $\beta\in (0,1)$ and (\ref{eqlokatzisPsi}), we arrive at
\[
Q(x)=Q(0)+o(1)+\mathcal{O}(\omega^{1-\beta})=o(1),\ \textrm{uniformly on}\ [0,\omega^{-\beta}],\ \textrm{as}\ n\to \infty.
\]
In light of the boundedness of $V_1'$ and (\ref{eqlokatzisPsi}), the above relation clearly implies that
\begin{equation}\label{eqagni10}
\varphi_1(x)=\lambda V_1'(x)+o(1),\ \textrm{uniformly on}\ [0,\omega^{-\beta}],\ \textrm{as}\ n\to \infty.
\end{equation}
Analogously, we can show that
\begin{equation}\label{eqagni11}
\varphi_2(x)=\lambda V_2'(x)+o(1),\ \textrm{uniformly on}\ [-\omega^{-\beta},0],\ \textrm{as}\ n\to \infty.
\end{equation}

\underline{Exchange of information between the inner and outer zones.}

Evaluating (\ref{eqtildenuSinhminus}) and (\ref{eqagni11}) at $x=-\omega^{-\beta}$, keeping in mind (\ref{eqasym}), we find
\begin{equation}\label{eqexchange1}
\nu_n \to -A\lambda.
\end{equation}
Similarly, from (\ref{eqtildemuSinhplus}) and (\ref{eqagni10}) for $x=\omega^{-\beta}$ we obtain
\begin{equation}\label{eqexchange2}
\mu_n\to A\lambda.
\end{equation}

\underline{Connecting the outer zones, bypassing the inner one.}

 Roughly speaking, so far we have shown that $(\varphi_1,\varphi_2)$ resembles a constant multiple of the counterexample from Proposition \ref{procounter} for large $n$. However, we have not yet made full use of the second limit in (\ref{eqproofContra}) (recall that the second estimate in (\ref{eqagni6}) was more than what was needed for (\ref{eqagni7})).

As in \cite[Prop. 2.2]{aftalionSour}, we wish to make a connection between the outer profiles (\ref{eqtildemuSinhplus}) and (\ref{eqtildenuSinhminus}).  To this end, testing (\ref{eqproofDirThm}) by $(V_1',V_2')$ over $(-\omega^{-\beta},\omega^{-\beta})$ yields
\begin{equation}\label{eqagni100}\begin{split}
\sum_{i=1}^{2}\left(\varphi_i'V_i'(-\omega^{-\beta})-\varphi_iV_i''(-\omega^{-\beta})\right)-\sum_{i=1}^{2}\left(\varphi_i'V_i'(\omega^{-\beta})-\varphi_iV_i''(\omega^{-\beta})\right) &    \\
  +\omega^2\sum_{i=1}^{2}\int_{-\omega^{-\beta}}^{\omega^{-\beta}}\varphi_iV_i'dx    =                                 \sum_{i=1}^{2} \int_{-\omega^{-\beta}}^{\omega^{-\beta}}g_iV_i'dx.&           \end{split}
\end{equation}
Recalling (\ref{eqasym}), (\ref{eqproofContra}), (\ref{eqgradg}), and making use of the estimates
\[
\omega^2\int_{-\omega^{-\beta}}^{\omega^{-\beta}}\varphi_iV_i'dx=\mathcal{O}(\omega^{2-\beta}),\  \int_{-\omega^{-\beta}}^{\omega^{-\beta}}g_iV_i'dx=o(\omega)\ (\textrm{cf.}\ (\ref{eqagni6})),
\]
keeping in mind that $\beta\in (0,1)$, it follows readily from (\ref{eqagni100}) that
\[
-A\varphi_2'(-\omega^{-\beta})-A\varphi_1'(\omega^{-\beta})=o(\omega)\ \textrm{as}\ n\to \infty.
\]
Hence, via (\ref{futureprecise}), we infer that
\begin{equation}\label{eqagni101}
  -\nu_n+\mu_n\to 0.
\end{equation}

\underline{Putting everything together and completion of the proof.}

By combining (\ref{eqexchange1}), (\ref{eqexchange2}) and (\ref{eqagni101}), we deduce that
\[
\lambda=0 \ \textrm{and}\ \mu_n, \nu_n\to 0\ \textrm{as}\ n \to \infty.
\]
Thus, by virtue of (\ref{eqinnerloc1}), (\ref{eqinnerMainLokatzis}), (\ref{eqtildemuSinhplus}), (\ref{eqtildenuSinhminus}), (\ref{eqagni10}) and (\ref{eqagni11}), we conclude that
\begin{equation}\label{eqlen111}
\|\varphi_1\|_{C_0(0,R_n)}\to 0,\ \varphi_1\to 0\ \textrm{in}\ C_{loc}(-\infty,0),\   \|\varphi_2\|_{C_0(-R_n,0)}\to 0,\ \varphi_2\to0 \ \textrm{in}\ C_{loc}(0,\infty).
\end{equation}
Therefore, to reach the desired contradiction with the first relation in (\ref{eqproofContra}),
it remains to establish the uniform vanishing of $\varphi_1$ on $[-R_n,0]$ and of $\varphi_2$ on $[0,R_n]$ as $n\to \infty$.
We will only deal with the case of $\varphi_2$ since that of $\varphi_1$ can be treated completely analogously. This task will be carried out in the remainder of the proof.

In light of (\ref{eqasym}), (\ref{eqproofDirThm}), (\ref{eqproofContra}) and (\ref{eqlen111}), we have
\[
-\varphi_2''+(V_1^2+\omega^2)\varphi_2=o(1)e^{-\theta x}, \ \textrm{uniformly on}\ [0,R_n],\ \textrm{as}\ n\to \infty;\ \varphi_2(R_n)=0.
\]
Now, let $K>0$ be such that
\[
V_1^2(K)=2\theta^2.
\]
Then, recalling that $V_1'>0$, by a standard barrier argument we obtain
\[
|\varphi_2(x)|\leq \left(|\varphi_2(K)|+o(1)\right)e^{-\theta (x-K)},\ x\in [K,R_n].
\]
Consequently, the desired uniform vanishing of $\varphi_2$ as $n\to \infty$ follows readily from (\ref{eqlen111}) and the above uniform bound.
\end{proof}

\begin{rem}\label{remorth}
  In analogy to \cite[Prop. 2.2]{aftalionSour} which dealt with the case $\omega=0$, if in Theorem \ref{thm} we make the further assumption
  \begin{equation}\label{eqorthrem}
      \int_{-R}^{R}(V_1'g_1+V_2'g_2)dx=0,
  \end{equation}
  then the stronger assertion
  \[
  \|\varphi\|_{C^0\left((-R,R) \right)}\leq C\|g\|_{C^0_\theta\left((-R,R) \right)}
  \]
holds. Indeed, the proof proceeds verbatim, the only essential difference being on how to estimate the righthand side of the corresponding relation to (\ref{eqagni100}).
To this end, we note that thanks to (\ref{eqorthrem}) and
\[
\|g_n\|_{C^0_\theta\left((-R_n,R_n) \right)}\to 0\ (\textrm{cf.}\ (\ref{eqproofContra})),
\]
 we would have
\[
\sum_{i=1}^{2}\int_{-\omega^{-\beta}}^{\omega^{-\beta}}g_iV_i'dx=-\sum_{i=1}^{2}\int_{\omega^{-\beta}<|x|<R}^{}g_iV_i'dx=\mathcal{O}\left(e^{-\frac{c}{\omega^\beta}}\right),
\]
which is more than enough to proceed.
\end{rem}

\end{document}